\documentclass[12pt]{amsart}
\usepackage{amsmath, color} \usepackage{amssymb,verbatim} \usepackage{amsthm}
\usepackage{amscd}
\usepackage{thmtools}
\usepackage{cleveref,fancyhdr,url}
\usepackage{tikz}
\usepackage[T1,T5]{fontenc}

\usepackage{tabularx}
\usepackage{tikz}
\usepackage{array}
\usepackage{tkz-graph}
\usepackage{amssymb,tikz}
\usepackage{tikz-cd}
\tikzcdset{scale cd/.style={every label/.append style={scale=#1},
		cells={nodes={scale=#1}}}}

\tikzset{smalltext/.style={"\textup{\small #1}" description}}
\usetikzlibrary{positioning}
\usepackage{tkz-euclide,amsmath} 
\usetikzlibrary{arrows.meta, bending}

\tikzset{hidden/.style = {thick, dashed}}

\def\c{{\mathfrak{c}}}
\def\F{{\mathbb{F}}}

\DeclareMathOperator{\st}{st}
\DeclareMathOperator{\reg}{reg}

\newcommand{\G}{\mathcal{G}}
\DeclareMathOperator{\dist}{d}
\DeclareMathOperator{\ngb}{N}

\theoremstyle{plain}
\newtheorem{theorem}{Theorem}[section]
\newtheorem{corollary}[theorem]{Corollary}
\newtheorem{proposition}[theorem]{Proposition}

\theoremstyle{definition}
\newtheorem{definition}[theorem]{Definition}

\newtheorem{lemma}[theorem]{Lemma}
\newtheorem{question}[theorem]{Question}
\newtheorem{example}[theorem]{Example}
\newtheorem{remark}[theorem]{Remark}

\def\fu{\mathbf{u}}
\def\f0{\mathbf{0}}

\def\fw{\mathbf{w}}
\def\fv{\mathbf{v}}

\def\fq{\mathbf{q}}

\def \G{\mathcal G}

\def \M{\mathcal M}

\def \F{\mathcal F}

\def \k{\Bbbk}
\def \D{\Delta}

\def \s{\Bbbk[x_1,x_2,\ldots,x_n]}
\def \la{\langle}
\def \ra{\rangle}
\newcommand{\dis}{\displaystyle}

\DeclareMathOperator{\supp}{{Supp}}
\DeclareMathOperator{\cp}{{Comp}}

\DeclareMathOperator{\maxd}{{maxdegree}}

\begin{document} 
\vspace*{-1cm}
\title{Good-vertex splittability and constraints for regularity}
\date{}

\author[H. Đào]{ Ho\`ai \DJ\`ao} 
\address[H. Đào]{Department of Mathematics, Virginia Tech, Blacksburg, VA 24061, USA}
\email{hoaidao@vt.edu}

\author[R. Nair]{Ritika Nair}
\address[R. Nair]{Department of Mathematical Sciences, Oklahoma State University, Stillwater, OK 74078, USA}
\email{ritika.nair@okstate.edu}
\urladdr{\url{https://sites.google.com/view/ritika-nair}}

\begin{abstract}

In this paper, we introduce the concept of a good vertex for a monomial ideal, which is a variable in the support of the ideal that satisfies a specific exchange property. A good vertex for a squarefree monomial ideal is a shedding vertex of the Alexander dual simplicial complex. We explore several nice consequences of having a good vertex, particularly for edge ideals of graphs. If almost all vertices are good, the edge ideal has a linear resolution. We also develop an inductive technique in terms of good vertex to show that several families of monomial ideals have linear resolution. This technique gives the stronger property of good-vertex splittability, which extends the idea of vertex splittability first defined by Moradi and Khosh-Ahang. We use some techniques developed in this paper to show that squarefree monomial ideals whose `complements' are small enough are good-vertex splittable. Lastly, we present some sufficient combinatorial conditions for edge ideals having regularity at most $3$ in terms of forbidden induced subgraphs of size $5$ or $6$, complementing the recent work by Kanno.

\end{abstract}

\maketitle

\vspace{-1cm}
\section{Introduction}

The Castelnuovo-Mumford regularity (or simply, regularity) of ideals in a polynomial ring is an important homological invariant. 
It is a well-known fact that a monomial ideal equigenerated in degree $d$ has a linear resolution if and only if the regularity of the ideal is $d$. We introduce the concept of a \textit{good vertex} which is a variable in the support of a monomial ideal that satisfies a specific exchange property (see \Cref{good_v}), and present the following inductive criterion that helps show that several families of ideals have linear resolution without explicitly computing the resolution:

{\renewcommand{\thetheorem}{\ref{gvs_linres}}
\begin{theorem}
  Let $\F$ be a family of equigenerated monomial ideals. Suppose every ideal $I \in \F$ is either a principal ideal or there is a good vertex $x \in \supp(I)$ such that both $(I:x)$ and $I_{-x}$ are in the same family $\F$. Then, any ideal in $\F$ has linear resolution. Here, $I_{-x}$ is the ideal generated by all the monomial generators of $I$ which do not contain $x$.
\end{theorem}
\addtocounter{theorem}{-1}
}

A good vertex is highly versatile. For a monomial ideal $I$ generated in degree $d$, a variable $x$ being good for $I$ is characterized by the fact that $(I:x)$ is generated in degree $d-1$ (see \Cref{prop_good}). For a squarefree monomial ideal, the good vertex naturally arises as the shedding vertex of the Alexander dual simplicial complex associated with the monomial ideal. In case of the edge ideal of a graph, a good vertex restricts the structure of the graph (\Cref{pro:distance,prop_good_degree}). In particular, such a graph is connected and the distance between any two vertices is at most $4$.\\

\Cref{gvs_linres} is stated for equigenerated ideals, and the proof can be found on Page 12. While the present paper is motivated by this result, it turns out that for any monomial ideal, the criteria in \Cref{gvs_linres} yield a stronger property that we call \emph{good-vertex splittability} (see \Cref{defn_gvs_family}), as it is an extension of the notion of vertex splittability introduced by Somayeh Moradi and Fahimeh Khosh-Ahang in \cite{MKA16}. The definition of good-vertex splittability has a subtle distinction and hence is a larger class of ideals which coincides with vertex splittability when the ideals are squarefree. Further, this larger class of ideals also has linear quotients and gives a Betti splitting, just as vertex splittable ideals. \\

In many situations in commutative algebra (respectively, combinatorics), one is interested in how various properties of a monomial ideal (respectively, simplicial complex) are preserved as one removes generators (respectively, faces) (see \cite{TY06},\cite{DNSVTW24}). A squarefree monomial ideal $I$ generated in degree $d$ is good-vertex splittable whenever its complement (the monomials of degree $d$ not in $I$) is small in a certain sense. We show this in two cases. \Cref{dual_vd} is the strengthening of a result by Naoki Terai and Ken-ichi Yoshida in \cite{TY06}. The next result (see \Cref{complement2}) shows that the ideal is good-vertex splittable when the complements on certain subsets of the support are small.\\

There is a characterization for edge ideals having regularity $2$ over any base field given by Fr\"oberg in terms of the graphs being co-chordal \cite{F90}. We show that if the edge ideal has too many good vertices, then the regularity is $2$: 

{\renewcommand{\thetheorem}{\ref{gvreg}}
\begin{theorem}
    If the edge ideal of a graph has at most $3$ bad vertices, it has a linear resolution. 
\end{theorem}
\addtocounter{theorem}{-1}
}

Attempts have been made to get characterization for edge ideals with higher regularities in terms of the structure of the graph. In \cite{FRG14}, Fern\'andez-Ramos and Philippe Gimenez gave a full characterization for edge ideals of bipartite graphs with regularity $3$. In \cite{Ne11}, Eran Nevo showed that a  graph which is gap-free and claw-free has regularity at most $3$. In \cite{Ba15}, Arindam Banerjee extended this result to show that edge ideal of a graph which is gap-free and cricket-free has regularity at most $3$. We present criteria for edge ideals having regularity at most $3$ in terms of forbidden induced subgraphs on $5$ (see \Cref{reg3graphfive}) or $6$ (see \Cref{reg3-graph-six}) vertices:

{\renewcommand{\thetheorem}{\ref{reg3graphfive}}
\begin{theorem}
   Let $G$ be a simple graph without isolated vertices such that the induced subgraph of $G$ on any five vertices  is neither $K_1 \sqcup 2K_2$ nor  $K_1 \sqcup P_4$.  Then the edge ideal of $G$ has regularity at most $3$. 
\end{theorem}
\addtocounter{theorem}{-1}
}

See \Cref{tab:six} for the forbidden subgraphs on $6$ vertices. We also comment on the connection of our work with the  recent work by Akane Kanno in \cite{Ka26} that gives topological characterization for edge ideals having regularity at most $3$.\\

The paper is organized as follows. Section 2 gives the background and sets some notations. The objective of this paper is threefold. We first introduce a good vertex and explore its consequences in Section 3. This also includes a nice combinatorial visualization for a good vertex for a squarefree monomial ideal. Second, Sections 4 and 5 elaborate on good-vertex splittability. Section 4 outlines the fact that these ideals have linear quotients adapting the proof similar to \cite{MKA16}. We include an independent proof of \Cref{gvs_linres} in this section that does not rely on linear quotients, which led to the development of this paper, and also give an alternative approach to \cite{MKA16} to show the property of Betti splitting. Section 5 presents some classes of good-vertex splittable ideals by demonstrating applications of the combinatorial visualization technique for a good vertex introduced in Section 3. Lastly, Section 6 has the criteria for edge ideals to have regularity at most $3$. Section 7 concludes with a few questions for further exploration.\\

 \textbf{Acknowledgements}: The authors are thankful to Hailong Dao for his valuable comments. The idea for this project developed when the second author was working on \cite{DNSVTW24} and is grateful to Anton Dochtermann, Jay Schweig, Adam Van Tuyl and Russ Woodroofe for their inputs. Some of the results in this paper were included in the Ph.D. dissertation thesis of Ritika Nair \cite{Na24} in 2024. 

\section{Preliminaries and Notations}

Throughout this paper, we assume $S=\k[x_1, \ldots, x_n]$ is a standard graded polynomial ring over a field $\k$ of characteristic zero and $[n]:=\{x_1,x_2,\ldots,x_n\}$. For a variable $x$, we set the notation $\k[\hat{x}]$ to refer to the polynomial ring $\k[x_i~|~x_i \ne x]$. 
\subsection{Monomial ideals}
For a monomial ideal $I \subset S$, we denote the unique minimal generating set of $I$ by $\G(I)$. When we say that $I$ is generated in degree $d$, we mean that $I$ is an \emph{equigenerated} monomial ideal, i.e., every monomial in $\G(I)$ has degree $d$. A monomial ideal is \emph{principal} if it is generated by a single monomial. We consider the zero ideal $\la 0 \ra$ as a monomial ideal. A monomial $\fu$ is \emph{squarefree} if $x_i^2$ does not divide $\fu$ for any $i$ and the monomial ideal $I$ is squarefree if the elements of $\G(I)$ are squarefree. Let $\M^{sq}_d$ denote the collection of all squarefree monomials of degree $d$ in $\s$ and $\M^{sq}_I$ denote the collection of all squarefree monomials in $I$.

The \emph{support} of a monomial $\fu$ is defined as $\supp(\fu):= \{x_i~|~x_i \text{ divides } \fu \}$ and the support of a monomial ideal $I$, $\supp(I)$, is the union of the supports of its minimal generators. We define $\maxd(I) := \max\{\deg(\fu)~|~\fu \in \G(I)\}$ where $\deg(\fu)$ is the degree of the monomial $\fu$. For a variable $x \in \supp(I)$, the ideal $I_{-x}:=\la \fu \in \G(I)~|~ x \notin \supp(\fu) \ra$.  \\

\begin{definition} \label{def_comp}
    Let $I \subset \s$ be a squarefree monomial ideal generated in degree $d$. 
    The \emph{complement ideal} (or simply, the complement) of $I$ is defined as
    \[\cp(I):= \la \fu \in \mathcal{M}^{sq}_{d} ~\setminus~ \G(I)~|~ \supp(\fu) \subset \supp(I) \ra \]
    If the collection of monomials in the generating set above is empty, then we set $\cp(I)=\la 0 \ra$.
    \end{definition}

The notion of a complementary ideal was defined by Keller Vandebogert in \cite{Va22}. Our definition is a special case due to the restriction placed on the support of the monomials generating $\cp(I)$.

\begin{example} \label{ex:basic}
    Let $I = \la x_1x_2, x_2x_3, x_3x_4 \ra \subset \k[x_1,x_2,x_3,x_4,x_5]$. Then, $\supp(I)=\{x_1,x_2,x_3,x_4\}$, $\cp(I) = \la x_1x_3,x_1x_4,x_2x_4 \ra$ and $I_{-x_1}= \la x_2x_3,x_3x_4 \ra$. 
    \end{example}

Given the minimal free resolution of a monomial ideal $I$,
\[
0 \xrightarrow{} \bigoplus_j S(-a_{p,j})^{\beta_{p,j}} \xrightarrow{} \cdots \xrightarrow{} \bigoplus_j S(-a_{1,j})^{\beta_{1,j}} \xrightarrow{} \bigoplus_j S(-a_{0,j})^{\beta_{0,j}} \xrightarrow{}I \xrightarrow{}0
\] 

\noindent  The numbers $\beta_{i,j}$ are called the graded \emph{Betti numbers} of $I$. The \emph{Castelnuovo-Mumford regularity} of $I$, $\displaystyle \reg(I) = \max_i\{j-i~|~ \beta_{i,j} \ne 0\}$. If $I$ is generated in degree $d$ and $\reg(I)=d$, then $I$ has a \emph{linear} resolution. We recall a simple observation below.

\begin{lemma} \label{lem_prodlinear}
    Let $I\subset S$ have a linear resolution. Then, $xI$ also has a linear resolution for any variable $x$.
\end{lemma}

The monomial ideal $I$ is said to have \emph{linear quotients} if there is an ordering $\fu_1<\fu_2< \cdots < \fu_t$ on the elements of $\G(I)$ such that $\dis \la \fu_1,\ldots, \fu_i\ra:\la \fu_{i+1}\ra$ is generated by variables for all $1 \le i \le t-1$. 
 
\begin{theorem}\cite[Lemma 5.2]{Fa04} \label{linquo_linres}
Let $I$ be an equigenerated monomial ideal of \\
$\s$. If $I$ has linear quotients, then $I$ has a linear resolution.    
\end{theorem}

In \cite{MKA16}, Somayeh Moradi and Fahimeh Khosh-Ahang defined the notion of \emph{vertex splittability} as follows: 

\begin{definition} \label{moradi - vs}
A monomial ideal $I \subset S=\s$  is called \emph{vertex splittable} if it can be obtained by the following recursive procedure.
\begin{enumerate}
    \item If $\fu$ is a monomial and $I = \la \fu \ra$ or $I=\la 0 \ra$ or $I = S$, then $I$ is a vertex splittable ideal.
    \item If there is a variable $x$ and vertex splittable ideals $I_1$ and $I_2$ of $\k[\hat{x}]$ so that $I=x I_1 + I_2$, $I_2 \subset I_1$, and $\G(I) = \G(xI_1)\sqcup \G(I_2)$, then $I$ is vertex splittable. The union here is disjoint and $x$ is called a \emph{splitting vertex}. 
\end{enumerate}
\end{definition}

\begin{theorem}\cite[Theorem 2.4]{MKA16} \label{vs_linquo}
    A vertex splittable ideal has linear quotients.
\end{theorem}

Another interesting consequence of vertex splittability is that of Betti splitting. For monomial ideals $I,J$ and $K$, $I = J + K$ is a \emph{Betti
splitting of $I$} if  the set of minimal monomial generators of $I$ is the disjoint union of the sets of minimal monomial generators of $J$ and $K$, and 
\begin{align*}
    \beta_{i,j} (I) = \beta_{i,j} (J) + \beta_{i,j} (K) + \beta_{i-1,j} (J \cap K),
\end{align*}
 for all $i, j \in \mathbb{N}$.

 \begin{theorem}\label{vs_betti}
     \cite[Theorem 2.8]{MKA16}
     If $I=xI_1+I_2$ is a vertex splitting for the monomial ideal $I$, then $I=xI_1+I_2$ is also a Betti splitting.
 \end{theorem}

\subsection{Simplicial complexes} Given the set $[n]=\{x_1,\ldots, x_n\}$ of vertices, an abstract \emph{simplicial complex} (or simply, a complex) $\Delta$ is a collection of subsets of $[n]$ which is closed under inclusion, i.e., $F \in \Delta$ and $E \subset F$ implies $E \in \Delta$. We assume $\{x_i\} \in \D$ for all $i$. Elements of $\Delta$ are called \emph{faces}. The maximal elements of $\Delta$ with respect to inclusion are called \emph{facets}. The \emph{dimension} of a face $F$, $\dim(F) = |F|-1$, where $|F|$ is the number of vertices in $F$ and a $d$-dimensional face is simply called a $d$-face. The dimension of $\Delta$, $\dim(\Delta)=\max\{|F|-1~|~ F \text{ is a facet of } \Delta\}$. $\Delta$ is said to be \emph{pure} if all its facets have the same  dimension. The \emph{$d$-skeleton} of a simplicial complex $\D$, $\D_d:=\{F \in \D~|~\dim(F) \le d\}$. A \emph{$k$-simplex} is a simplicial complex having a single facet of dimension $k$. Since a simplicial complex is closed under inclusion, it is sufficient to mention the facets to describe the complex. We often denote a face $\{x_{i_1},x_{i_2},\ldots,x_{i_t}\}$ simply as $x_{i_1}x_{i_2}\cdots x_{i_t}$. \\

The vertex $x$ is a \emph{shedding vertex} for $\D$ if for every facet $F$ containing $x$, there is a vertex $y \notin F$ such that $(F \setminus \{x\}) \cup \{y\}$ is still a face of $\D$.\\

Given a squarefree monomial ideal $I\subset \s$, one can associate a complex to $I$ via Stanley-Reisner theory. We require the notion of the Alexander dual complex of $I$, $\D_I^\vee$, whose facets are the complements of the monomials in $\G(I)$, i.e., $\D_I^\vee:=\left\la  \dfrac{x_1x_2\cdots x_n}{\fu}~|~\fu \in \G(I) \right\ra$. This correspondence between squarefree monomial ideals in $\s$ and simplicial complexes on $[n]$ is a bijection.\\

 The motivation for the definition of vertex splittability \cite{MKA16} by Moradi and Khosh-Ahang stems from the need to analyze the consequence of the important combinatorial concept of vertex decomposability of a simplicial complex on the Alexander dual ideal. Vertex decomposability was first introduced by J. Scott Provan and Louis J.Billera \cite{PB80} in 1980 for pure complexes and later extended to non-pure complexes by Anders Bj\"orner and Michelle L. Wachs \cite{BW97}. We skip the definition here, and state the relevant connection below.

\begin{theorem}\cite[Theorem 2.3]{MKA16} A squarefree monomial ideal $I$ is vertex splittable if and only if $\D_I^\vee$ is vertex decomposable.
\end{theorem}

\subsection{Graph theory}
A finite simple graph $G$ on the vertex set $[n]$ is a $1$-dimensional simplicial complex. The set of $1$-faces is the edge set $E(G)$. The \emph{edge ideal} of $G$, $I(G):= \la x_ix_j~|~ \{x_i,x_j\} \in E(G) \ra \subset \s$. An \emph{induced subgraph} $H$ of $G$ is the graph on a subset $W\subset [n]$ with the edge set $E(H) = \{\{x_i,x_j\} \in E(G)~|~x_i,x_j \in W\}$. $G \setminus H$ is the induced subgraph on $[n]\setminus W$.\\

The \emph{complement} of the graph $G$, $G^{\c}$ is defined on the same vertex set $[n]$ and has the edge set $E(G^{\c}):=\{ \{x_i,x_j\}~|~ i \ne j, \{x_i,x_j\} \notin E(G) \}$. Further, $(G^{\c})^{\c}=G$ and for any induced subgraph $H$ of $G$, $H^{\c}$ is an induced subgraph of $G^{\c}$.\\

 Next, we recall the definitions of specific types of graphs. A \emph{path graph} $P_n$ is a graph on $[n]$ having edge set $E(P_n)=\{\{x_i,x_{i+1}\}~|~ 1 \le i \le n-1\}$, and the \emph{$n$-cycle} $C_n$ has edge set $E(P_n) \cup \{x_1,x_n\}$. A graph is called \emph{chordal} if every cycle has an induced subcycle of length $3$. Every induced subgraph of a chordal graph is chordal. A graph $G$ on $4$ vertices is called a \emph{gap} if $G=K_2\sqcup K_2$, where $K_n$ is the complete graph on $n$ vertices. A graph $G$ is said to be \emph{gap-free} if it does not have an induced subgraph which is a gap.\\

Ralf Fr\"{o}berg gave the following celebrated result that provides a combinatorial characterization for squarefree monomial ideals generated in degree $2$ having a linear resolution.

\begin{theorem}  \cite[Theorem 1]{F90} \label{thm_chordal}
    A squarefree monomial ideal generated in degree $2$ has a linear resolution if and only if it is the edge ideal of a graph whose complement is chordal. 
\end{theorem} 

Moradi and Khosh-Ahang strengthened this result.

\begin{theorem} \cite[Theorem 3.6, Corollary 3.8]{MKA16} \label{thm_chordal_vs}
    A squarefree monomial ideal generated in degree $2$ is vertex splittable if and only if it is the edge ideal of a graph whose complement is chordal. 
\end{theorem}

For two vertices $x_i,x_j$, the distance $\dist_G(x_i,x_j)$ is the number of edges of a shortest path from $x_i$ to $x_j$. $G$ is said to be connected if for any pair of vertices, there is a path between them. When $\dist_G(x_i,x_j)=1$, $x_j$ is a \emph{neighbor} of $x_i$. The set of all neighbors of $x_i$ in $G$ is denoted by $\ngb_G(x_i)$ and the \emph{degree} of a vertex $x_i$, $\deg(x_i) = |\ngb_G(x_i)|$. A vertex $x_i$ is isolated if  $\deg(x_i)=0$. The collection of vertices $\{x_i\}\cup \ngb_G(x_i)$ is denoted by $\st_G(x_i)$, and called the \emph{star} of $x_i$ in $G$. We denote the subgraph of $G$ on $[n]\setminus \st_G(x_i)$ by $G \setminus \st_G(x_i)$.\\

When $G$ does not have isolated vertices, $\cp(I(G))=I(G^\c)$. Since we use $x_i$'s to refer to the vertices of a complex on $[n]$ as well as the variables of $\s$, henceforth we refer to the $x_i$'s as variables or vertices interchangeably. 

\section{Good vertex\nopunct}

\begin{definition}\label{good_v}
     Let $I \subset  \s$ be a monomial ideal, $x \in \supp(I)$ and $\fu \in I$ be a monomial. The pair $(x,\fu)$ is said to be a \emph{good pair} if $x$ satisfies the ``exchange'' condition that there exists some $y \in \supp(\fu)$ such that $\dis \frac{\fu}{y} \cdot x \in I$. Otherwise, $(x,\fu)$ is a \emph{bad pair}. If $x$ satisfies the exchange condition for every monomial $\fu \in I$, then $x$ is called a \emph{good vertex} for $I$. Else, it is a bad vertex. \\
Note that for any $x \in \supp(I)$, to check if $x$ is a good vertex for $I$, it is sufficient to check the exchange condition on monomials $\fu \in \G(I_{-x})$. In particular, if $x \in \supp(\fu)$ for every $\fu \in \G(I)$, then $x$ is good for $I$.

\end{definition}

\begin{example}\label{eg1_good}
     Let $I = \la ae, b^2cd, de \ra \subset \k[a,b,c,d,e,f]$.
    Here $e$ is a good vertex for $I$ since $ \dis I_{-e}=\la b^2cd \ra$ and  $ \dis \frac{b^2cd}{b} \cdot e = bcde \in I$. On the other hand, for the vertex $a$, $I_{-a}=\la b^2cd, de \ra$ and $(a,b^2cd)$ is a bad pair.
\end{example}

\begin{proposition}\label{prop_good} 
Let $I$ be generated in degree $d$. For $x \in \supp(I)$ and a monomial $\fu \in \G(I)$, $(x,\fu)$ is a bad pair if and only if $\fu$ is a minimal generator of $(I:x)$. \\
Thus, $x$ is good for $I$ if and only if $(I:x)$ is generated in degree $d-1$ only. 
\end{proposition}

\begin{proof}
    Let $x \in \supp(I)$ and $\fu \in \G(I)$ for a monomial ideal $I$ generated in degree $d$. We have $\fu \in I \subset (I:x)$. Without loss of generality, let $x \notin \supp(\fu)$. If $\fu$ is a minimal generator of $(I:x)$, then for any $y \in \supp(\fu)$, $\dis \frac{\fu}{y} \cdot x \notin I$, i.e., $(x,\fu)$ is a bad pair.\\ 
    Next, suppose $\fu$ is not a minimal generator of $(I:x)$. Then there exists a minimal generator $\fv$ of $(I:x)$ with $\deg(\fv) < d$ such that $\fv$ divides $\fu$. Let $\fu = \fv \cdot \fw$, $\deg(\fw) \ge 1$. For any $ y \in \supp(\fw)$ (note that $y \ne x$), $\dis \frac{\fu}{y} \cdot x =  \frac{(\fv \cdot\fw) }{y} \cdot x = (\fv \cdot x) \cdot \frac{\fw}{y} \in I$ as $\fv \cdot x \in I$. So, $(x, \fu)$ is a good pair.
\end{proof}

Given a monomial ideal $I$, the property of having a good vertex is not preserved upon taking the colon operation.

\begin{example} \label{ex_colonideal}
 Let $I=\la abd,abf,ace,adc,aef,bde,bcf,bce,cdf,def \ra$. Then any vertex is good for $I$. However, $J=(I:a)$ is the edge ideal of a $5$-cycle, so no vertex is good for $J$ as demonstrated in \Cref{ex_goodvertex}. 
\end{example}

\textbf{A combinatorial visualization}: Now we look at the special case where $I$ is a squarefree monomial ideal with $\supp(I)=[n]$. We give a geometric realization for the concept of a bad pair by coloring the faces of a simplex. Let $\D$ be the simplex on $n$ vertices. For every $d \ge 0$, there is a bijection from the collection of all $d$-faces to the monomials in $\M^{sq}_{d+1}$. For each $d$, color the faces $x_{i_1}x_{i_2} \cdots x_{i_d}$ of $\D$ corresponding to degree $d$ squarefree monomials in $I$ in blue and the others in red. For $x \in \supp(I)$ and a monomial $\fu \in \M^{sq}_I$ of degree $d$ with $x \notin \supp(\fu)$, define $\D_{(x,\fu)} \subset \D$ as the $(d-1)$-skeleton of the subcomplex on the $d+1$ vertices $\supp(\fu) ~\bigcup~\{x\}$. Let us look at the  facets of $\D_{(x,\fu)}$. There are $d+1$ of them given by the $(d-1)$-faces of $\D_{(x,\fu)}$. The facet of $\D_{(x,\fu)}$ corresponding to $\fu$ is colored blue since $\fu \in I$. Every other facet corresponds to a monomial of the type $  \frac{\fu}{y}\cdot x$ for some variable $y \in \supp(\fu)$. Hence, $(x,\fu)$ is a bad pair for $I$ if and only if $\D_{(x,\fu)}$ has only one blue facet corresponding to $\fu$ ``opposite'' to the vertex $x$, and the remaining $d$ many facets are all red. \\

For an ideal generated in degree $d$, it is sufficient to consider the coloring of the $(d-1)$-skeleton of the simplex on $[n]$ vertices. For any $n \ge 2$, the $1$-skeleton of the simplex on $[n]$ has edges corresponding to monomials in $\M^{sq}_2$. So, when $I$ is the edge ideal of a graph $G$, we get the visualization by coloring the edges of $G$ in blue and those in $G^{\c}$ in red. For a vertex $x$ of $G$, we abuse language and say $x$ is a good (or bad) vertex for $G$ to refer to the same property for the edge ideal $I(G)$. The following are some examples to illustrate this realization in the case of edge ideals of $n$-cycles and their complement graphs.

\begin{example}  \label{ex_goodvertex} \begin{enumerate}
 \item  Let $I$ be the edge ideal of $C_n$, the $n$-cycle. It is easy to see that every vertex is good for $C_n$ when $n<5$. 
    The case of $I=I(C_5)=\la x_1x_2,x_2x_3,x_3x_4,x_4x_5,x_1x_5\ra \subset \k[x_1,x_2,x_3,x_4,x_5]$ is illustrated below. 
    \begin{center}

\tikzset{every picture/.style={line width=0.75pt}} %set default line width to 0.75pt        

\begin{tikzpicture}[x=0.65pt,y=0.65pt,yscale=-1,xscale=1]

%Shape: Polygon [id:dp1313847156861967] 
\draw  [color={rgb, 255:red, 74; green, 144; blue, 226 }  ,draw opacity=1 ][line width=1.5]  (338.83,100.93) -- (317.82,170.44) -- (247.7,171.16) -- (225.37,102.1) -- (281.69,58.69) -- cycle ;
%Straight Lines [id:da011509418691514783] 
\draw [color={rgb, 255:red, 208; green, 2; blue, 27 }  ,draw opacity=1 ][line width=1.5]    (281.69,58.69) -- (317.82,170.44) ;
%Straight Lines [id:da44829968436523004] 
\draw [color={rgb, 255:red, 208; green, 2; blue, 27 }  ,draw opacity=1 ][line width=1.5]    (225.37,102.1) -- (338.83,100.93) ;
%Straight Lines [id:da1783593136810857] 
\draw [color={rgb, 255:red, 208; green, 2; blue, 27 }  ,draw opacity=1 ][line width=1.5]    (281.69,58.69) -- (247.7,171.16) ;
%Straight Lines [id:da4483390067130697] 
\draw [color={rgb, 255:red, 208; green, 2; blue, 27 }  ,draw opacity=1 ][line width=1.5]    (225.37,102.1) -- (317.82,170.44) ;
%Straight Lines [id:da9864862728370414] 
\draw [color={rgb, 255:red, 208; green, 2; blue, 27 }  ,draw opacity=1 ][line width=1.5]    (338.83,100.93) -- (247.7,171.16) ;

% Text Node
\draw (275,45.24) node [anchor=north west][inner sep=0.75pt]    {$x_{1}$};
% Text Node
\draw (207.95,95.38) node [anchor=north west][inner sep=0.75pt]    {$x_{2}$};
% Text Node
\draw (242.13,174.56) node [anchor=north west][inner sep=0.75pt]    {$x_{3}$};
% Text Node
\draw (310.87,172.68) node [anchor=north west][inner sep=0.75pt]    {$x_{4}$};
% Text Node
\draw (341.05,94.49) node [anchor=north west][inner sep=0.75pt]    {$x_{5}$};

\end{tikzpicture}

\end{center}
Consider the pair $(x_1,~x_3x_4)$.
  $\D_{(x_1,~x_3x_4)}$ is the subcomplex on $\{x_1,x_3,x_4\}$. Among the facets of $\D_{(x_1,~x_3x_4)}$, the blue edge corresponds to $x_3x_4$ opposite to $x_1$ and the remaining two edges are red as $x_1x_3$ and $x_1x_4$ are both not in $I$. So, this is a bad pair for $I$. In fact, we can extend this observation for any $n \ge 5$. $I= \la \{x_ix_{i+1}~|~1 \le i \le n-1\},~x_1x_n\ra$. For any vertex $x_i$, the pair $(x_i, x_{i+2}x_{i+3})$ (where we suitably take the index modulo $n$) is a bad pair for $I$. \\

    \item  For any $n \ge 4$, let $I$ be the edge ideal of $G=(C_n)^{\c}$, the graph complement of $C_n$. For any vertex $x_i$, we see that the pair $(x_i,x_{i-1}x_{i+1})$ is a bad pair for $I$. So, every vertex is a bad vertex for $I$. The case $n=4$ is shown below.

    \begin{center}

\tikzset{every picture/.style={line width=0.75pt}} %set default line width to 0.75pt        

\begin{tikzpicture}[x=0.65pt,y=0.65pt,yscale=-1,xscale=1]

%Shape: Diamond [id:dp4769877272569445] 
\draw  [color={rgb, 255:red, 208; green, 2; blue, 27 }  ,draw opacity=1 ][line width=1.5]  (297.18,81.82) -- (353.58,138.1) -- (297.18,194.38) -- (240.79,138.1) -- cycle ;
%Straight Lines [id:da6134844923162297] 
\draw [color={rgb, 255:red, 74; green, 144; blue, 226 }  ,draw opacity=1 ][line width=1.5]    (240.79,138.1) -- (353.58,138.1) ;
%Straight Lines [id:da736803878847786] 
\draw [color={rgb, 255:red, 74; green, 144; blue, 226 }  ,draw opacity=1 ][line width=1.5]    (297.18,81.82) -- (297.18,194.38) ;

% Text Node
\draw (289.73,68.5) node [anchor=north west][inner sep=0.75pt]    {$x_{1}$};
% Text Node
\draw (223.11,131.47) node [anchor=north west][inner sep=0.75pt]    {$x_{2}$};
% Text Node
\draw (290.73,194.3) node [anchor=north west][inner sep=0.75pt]    {$x_{3}$};
% Text Node
\draw (355.89,133.47) node [anchor=north west][inner sep=0.75pt]    {$x_{4}$};

\end{tikzpicture}

\end{center}
    
\end{enumerate}

\end{example}

\begin{proposition}\label{prop_commonred}
Given a squarefree monomial ideal $I$ generated in degree $d \ge 2$, let $(x,\fu)$ and $(y,\fv)$ be two bad pairs for $x \ne y \in [n]$, $\fu,\fv \in \G(I)$. Then, $\D_{(x,\fu)}$ and $\D_{(y,\fv)}$ can have at most one common red facet. 
\end{proposition}

\begin{proof} 
The pairs $(x,\fu)$ and $(y,\fv)$ are bad implies $x \notin \supp(\fu)$ and $y \notin \supp(\fv)$.  Suppose $\dis \D_{(x,\fu)}$ and $\dis \D_{(y,\fv)}$ have at least two distinct common red facets, say, corresponding to the monomials $\fq$ and $\textbf{p}$. We first consider $d>2$. We have  $\dis \frac{\fu}{x'}\cdot x = \fq = \frac{\fv}{y'}\cdot y$ and so, both $x$ and $y$ divide $\fq$. Let $\fq=\fw \cdot x \cdot y$. Then, $\fu=\fw\cdot y\cdot x'$ and $\fv = \fw \cdot x \cdot y'$. Note that $x' \ne y'$ else, $(x,\fu)$ and $(y,\fv)$ will be good pairs. Similarly from $\textbf{p}$, we get for some $\tilde{x}\ne y, ~\tilde{x} \ne x'$, and for some $\tilde{y}\ne x, ~\tilde{y} \ne y'$,  \[\dis \frac{\fu}{\tilde{x}}\cdot x = \frac{\fv}{\tilde{y}}\cdot y~ \text{ gives }~  \frac{\fw\cdot y\cdot x'}{\tilde{x}}\cdot x = \frac{\fw \cdot x \cdot y'}{\tilde{y}}\cdot y \]
     and hence, $\dis \tilde{y}\cdot x' = \tilde{x} \cdot y'$ resulting in either $x'=\tilde{x}$ or $x' = y'$, both contradicting our assumptions.\\
     If $d=2$, $\fq= x\cdot y$ gives $\fu = y\cdot x'$ and $\fv=x\cdot y'$ for $x' \ne y'$. Then, having another common red facet is possible only if \[\dis \frac{\fu}{y}\cdot x = \frac{\fv}{x}\cdot y\] giving $x\cdot x'=y \cdot y'$, i.e., $x=y', y=x'$ contradicting the squarefree property of $\G(I)$.
\end{proof}

\subsection{Good vertex and edge ideals} We now mention more consequences of having a good vertex for the edge ideal $I(G)$ of a graph $G$. 

The following proposition can be easily verified by the reader.

\begin{proposition} \label{prop_induced_bad}
    Let $I(G)$ be the edge ideal of a graph $G$ and $x \in \supp(I)$. Suppose $H$ is an induced subgraph of $G$ with $x \in \supp(I(H))$. If $x$ is a bad vertex for $H$, then $x$ is a bad vertex for $G$. The converse is not true. However, if $x$ is a bad vertex for $G$, there exists some induced subgraph $H$ of $G$ such that $x$ is bad for $H$. 
\end{proposition}

\begin{proposition} \label{pro:distance}
    Let $I=I(G)$ be the edge ideal of a graph on $[n]$ without any isolated vertex, and $x$ be a good vertex for $I$. Then,
    \begin{enumerate}
    \item  $G$ is connected.
        \item  For a vertex $y$ of $G$, $\dist_G(x,y) \leq 2$. 
        \item  $G \setminus \st_G(x)$ is either empty or is a union of isolated vertices. 
        \item  If $y,z$ are distinct vertices of $G$ and $y,z \in \dis [n] \setminus  \st_G(x)$, then $\dist_G(y,z) \leq 4$.
    \end{enumerate} 
\end{proposition}
\newpage
\begin{proof} \noindent 
\begin{enumerate}
\item  Suppose $G$ is not connected. Let $C_1$ and $C_2$ be any two connected components of $G$ and suppose $x \in C_1$. Then, for any edge $\{y,z\}$ in $C_2$, $(x,yz)$ is a bad pair for $I(G)$, giving a contradiction.
    \item   Let $y$ be a vertex of $G$ such that $y \notin \ngb_G(x)$. Then, for any vertex $z$ with $yz \in I(G)$, $xz$ is also in $I(G)$ as $x$ is a good vertex for $G$ and $xy \notin I(G)$. This gives $\dist_G(x,y) =2$.
    \item  Let $H=G \setminus \st_G(x)$. Suppose $H$ has at least two vertices, say $y$ and $z$. So, $xy \notin I(G)$ and $xz \notin I(G)$. Hence, $yz \notin I(G)$ as $x$ is a good vertex for $G$. So, $H$ is a collection of isolated vertices.
    \item  If $y$ and $z$ are vertices in $[n] \setminus \st_G(x)$, by the triangle inequality for distances in $G$ and by part (2), $\dist_G(y,z) \leq \dist_G(y,x) + \dist_G(x,z) = 4$.
\end{enumerate}
\end{proof}

\begin{remark}
    For a connected graph $G$, even if every pair of vertices $x,y$ satisfies $\dist_G(x,y) \leq 2$, $G$ may not have a good vertex. For example, $G=C_5$ (see \Cref{ex_goodvertex} (1)).
\end{remark}

\begin{proposition} \label{prop_good_degree}
   Suppose $x$ is a good vertex for a graph $G$. Let the maximum degree of any vertex of $G$ be $d$ and $W=\{y \in [n]~|~\deg(y)=d\}$. Then, one of the following is true:
   \begin{enumerate}
       \item $\deg(x)=d$ or
       \item  $\dis x \in \bigcap_{y ~\in~ W} \ngb_G(y)$, and $\dist_G(y,z) \le 2$ for every $y,z \in W$.
   \end{enumerate} 
\end{proposition}
\begin{proof}
   Let $x$ be a good vertex for $G$ with $\deg(x)<d$. Suppose there exists some vertex $y$ of maximum degree $d$ such that $x \notin \ngb_G(y)$. Then, $xw \in I(G)$ for all $w \in \ngb_G(y)$ as $x$ is a good vertex. This implies $\deg(x) \ge d$, giving a contradiction. So, $\dis x \in \bigcap_{y ~\in~ W} \ngb_G(y)$. As a consequence, for vertices $y\ne z$ in $W$, both $xy, xz \in I(G)$. If $y \notin \ngb_G(z)$, a shortest path from $y$ to $z$ is given by $\{y,x\},\{x,z\}$, giving $\dist_G(y,z) \le 2$.
\end{proof}

\subsection{A splitting vertex is a good vertex} For a squarefree monomial ideal $I$ with $x \in \supp(I)$, $\fu \in \G(I_{-x})$ if and only if $x$ is in the facet of $\D_I^\vee$ corresponding to $\dfrac{x_1x_2\cdots x_n}{\fu}$. 

  \begin{proposition}
       For a squarefree monomial ideal $I$, $x$ is a good vertex for $I$ if and only if $x$ is a shedding vertex for the Alexander dual complex of $I$. 
   \end{proposition}
   \begin{proof} Note that if $I_{-x}=\la0\ra$, then $x \in \supp(\fu)$ for all $\fu \in \G(I)$ and so, $x \notin F$ for any facet $F$ of $\D_I^\vee$. So, assume $I_{-x} \ne \la0\ra$.
       $x$ is a good vertex for $I$ if for every $\fu \in \G(I_{-x})$, $~\dis \frac{\fu}{y}\cdot x \in I$ for some $y \in \supp(\fu)$. This happens if and only if for every facet $F$ of $\D_I^\vee$ containing $x$, there is some $y \notin F$ such that $F \setminus \{x\} \cup \{y\} \in \D_I^\vee$.
   \end{proof}

This points towards the obvious connection between the notion of a good vertex and a splitting vertex.

\begin{proposition} \label{prop_splgood}
    A splitting vertex is a good vertex.
\end{proposition}

\begin{proof}
    Let $I\ne\la0 \ra$ be a monomial ideal with a splitting vertex $x$. Then, $I=xI_1 + I_2$, with $I_2 \subset I_1$ and $\G(I) = \G(xI_1)\sqcup \G(I_2)$. Note that $x \in \supp(I)$, as if $I_1 = \la 0 \ra$, then $I=\la0\ra$. Let $\fu \in I_2$, i.e., $x \nmid \fu$. We show that $(x,\fu)$ is a good pair. $\fu \in I_1$ implies $\fu=\fv \cdot\fw$ for some  $\fv \in I_1$ and monomial $\fw \in \k[\hat{x}]$ and since $\fu$ is a minimal generator of $I$, $\fu \ne \fv$ as $\fv \cdot x$ is a minimal generator of $I$. So, $\fw \ne 1$. For any $y \in \supp(\fw)$,  we get $ \dis \frac{\fu}{y} \cdot x = (\fv \cdot x)\left(\frac{\fw}{y}\right) \in I$, as required.    \end{proof}

 Note that a good vertex need not be a splitting vertex, unless it gives a splitting of the ideal.
      \begin{example}
          Let $I = \la x_1x_2,x_1x_3,x_1x_4,x_1x_5,x_2x_3,x_4x_5 \ra \subset \k[x_1,x_2,\ldots,x_5]$. Here, $x_1$ is a good vertex for $I$ but $x_1$ is not a splitting vertex.
      \end{example}

Hence, to check for vertex splittability of an ideal, we need to look for a good vertex giving a splitting.

\section{Good-vertex splittability \nopunct}

We are now ready to define the notion of good-vertex splittability.

\begin{definition}\label{defn_gvs_family}
    Let $\F$ be a family of monomial ideals such that every ideal $I \in \F$ is either a principal monomial ideal or there exists a good vertex $x \in \supp(I)$ such that both $(I:x)$ and $I_{-x}$ are in the same family $\F$. Then, $I=x\cdot(I:x)+I_{-x}$ is a good-vertex splitting and the ideals in $\F$ are said to be \it{good-vertex splittable}. In particular, for a good-vertex splitting $I=x\cdot(I:x)+I_{-x}$, both $(I:x)$ and $I_{-x}$ are good-vertex splittable. We set the convention that $I=S$ is good-vertex splittable.
\end{definition}

Note that when $I=x\cdot(I:x)+I_{-x}$ is a good-vertex splitting, \\$\G(I)= \G(x(I:x)) \sqcup \G(I_{-x})$, where $\G(x(I:x))=\{x\cdot \fu~|~\fu \in \G((I:x)) \}$. 

\begin{example}\label{gvs_ex} \begin{enumerate}
    \item An ideal generated by variables is good-vertex splittable, with any of the generating variables being the good vertex giving the splitting. For $\dis I = \la x_1,\ldots,x_m \ra \subset S=\s$, $I=x_1S +  \la x_2,\ldots, x_m\ra$ is a good-vertex splitting by induction.\\

   \item  Let $\dis I=\la x_1^2x_2,x_1x_2^2,x_1x_2x_3,x_1x_2x_4,x_1x_3x_4,x_2x_3x_4,x_3^2x_4,x_3x_4^2 \ra \subset S=\k[x_1,x_2,x_3,x_4]$. It is easy to see that $x_1$ is a good vertex for $I$, with $J=(I:x_1)=\la x_1x_2,x_2^2,x_2x_3,x_2x_4,x_3x_4\ra \subset S$, and $I_{-x_1}=\la x_2x_3x_4,x_3^2x_4,x_3x_4^2\ra \subset \k[x_2,x_3,x_4]$.\\
     $x_2$ is a good vertex for $J$ and $(J:x_2)=\la x_1,x_2,x_3,x_4 \ra$, $J_{-x_2}=\la x_3x_4\ra$, both of which are good-vertex splittable. \\
     Next, $x_2$ is a good vertex for $I_{-x_1}$. Let $L=I_{-x_1}$. $(L:x_2)=\la x_3x_4\ra$ is good-vertex splittable. $N=L_{-x_2}=\la x_3^2x_4,x_3x_4^2 \ra \subset \k[x_3,x_4]$ has good vertex $x_3$ with $N_{-{x_3}}=\la 0 \ra$ and $(N:x_3)=\la x_3x_4,x_4^2 \ra = x_3 \la x_4 \ra + \la x_4^2 \ra$, both of which are good-vertex splittable.

     \end{enumerate}
\end{example}

 Note that the definition of good-vertex splittability does not assume $\supp(I)=n$. So, if $I \subset S$ is good-vertex splittable, then $I$ is good-vertex splittable when considered as an ideal in $S[x_{n+1}]$.

\begin{lemma}\label{vs_implies_gvs}
    A vertex splittable ideal is good-vertex splittable.
\end{lemma}

\begin{proof}
    Let $\F$ be the family of vertex splittable ideals. Let $I \in \F$, $I \subseteq S$. If $I$ is a principal ideal or $I=S$, then $I$ is good-vertex splittable. Else, $I=xI_1+I_2$ where $I_1$ and $I_2$ are vertex splittable ideals in $\k[\hat{x}]$, $I_2 \subset I_1$, and $\G(I) = \G(xI_1)\sqcup \G(I_2)$, where $x$ is a splitting vertex. By \Cref{prop_splgood}, $x$ is a good vertex. Now, $\G(I_{-x})=\{ \fu \in \G(I)~|~ x \notin \supp(\fu) \} = \G(I_2)$. Next, we show $I_1=(I:x)$. A monomial $\fu \in (I:x)$ if and only if $x \fu \in I=x I_1+I_2$ if and only if $x \fu \in  x I_1$ as $I_2 \subset  \k[\hat{x}]$, and this is possible if and only if $\fu \in I_1$. Further by definition, both $I_1,~I_2 \in \F$. Hence, $I=xI_1+I_2$ is a good-vertex splitting.
      \end{proof}

Note that when $I=xI_1+I_2$ is a vertex splitting as per the definition in \cite{MKA16}, the degree of the splitting vertex $x$ in the generators of $I$ is at most $1$ as $I_1, I_2 \subseteq \k[\hat{x}]$. This need not be the case with good-vertex splitting as given a good-vertex splitting $I=x\cdot(I:x)+I_{-x}$ for a monomial ideal $I \subseteq \s$, $I_{-x} \subseteq \k[\hat{x}]$ by definition but $(I:x)$ can be in $\s$, as in \Cref{gvs_ex} (2) and \Cref{ex_borel}. Thus, the family of good-vertex splittable ideals is a larger class of ideals. When $I$ is squarefree, the two notions are equivalent.\\

The consequences of vertex splittability from \cite{MKA16} that are outlined in Section 2.2 apply for a good-vertex splittable ideal as well. 

\begin{theorem} \label{gvs_linquo}
    Let $I$ be good-vertex splittable, then $I$ has linear quotients.
\end{theorem}

\begin{proof}
    The proof technique in \cite[Theorem 2.4]{MKA16} applies here as well, by performing induction on the lexicographical order of the pair $(d,\ell)$, where $d:= \maxd(I)$ and $\ell:=| \supp(I)|$, so we skip the proof here.
\end{proof}

Corollary 2.5 in \cite{MKA16} states that an equigenerated vertex splittable ideal has a linear resolution and it is obtained by combining \Cref{linquo_linres} and \Cref{vs_linquo}. We provide an alternative direct approach to show that a good-vertex splittable ideal has a linear resolution, by applying the following Lemma by Hailong Dao, Craig Huneke and Jay Schweig  \cite[Lemma 2.10]{DHS13} giving bounds for the regularity of a monomial ideal.

\begin{lemma}\label{exact}
Let $I \subseteq S$ be a monomial ideal, and let $x \in \supp(I)$. Then
\[
\reg(I) \leq \max\{ \reg (I : x) + 1, \reg (I,x)\}. 
\]
Moreover, if $I$ is squarefree,  $\reg(I)$ is {\it equal} to one of these terms.
\end{lemma}

\begin{theorem} \label{gvs_linres} 
    Let $\F$ be a family of equigenerated monomial ideals. Suppose every ideal $I \in \F$ is either a principal ideal or there is a good vertex $x \in \supp(I)$ such that both $(I:x)$ and $I_{-x}$ are in the same family $\F$. Then, any ideal in $\F$ has linear resolution. Here, $I_{-x}$ is the ideal generated by all the monomial generators of $I$ which do not contain $x$.
\end{theorem} 

In short, an equigenerated good-vertex splittable ideal has a linear resolution.
    \begin{proof}
     Let $I$ be a good-vertex splittable ideal. We prove the theorem by induction on the lexicographical order of the pair $(d,\ell)$, where $d=\maxd(I)$ and $\ell$ is the cardinality of $\supp(I)$. The case $\ell=1$ is a principal ideal and has linear resolution. If $d=1$, then $I$ is generated by variables, and is resolved minimally by the Koszul complex which is linear. Assume that any good-vertex splittable ideal generated in degree $a$ and having support of size at most $b$ such that $(a,b)<(d,\ell)$ has linear resolution. Suppose $I$ is a good-vertex splittable ideal generated in degree $d$ with $|\supp(I)|=\ell$. Since $I$ is generated in degree $d$, $\reg(I)\ge d$, so it suffices to show that $\reg(I)\le d$. Let $x$ be the good vertex giving the splitting. $(I:x)$ is good-vertex splittable and is generated in degree $d-1$ by \Cref{prop_good}. $I_{-x}$ is good-vertex splittable and has $|\supp(I)|<\ell$. Suppose $I_{-x} \ne \la 0 \ra$. By induction, both $(I:x)$ and $I_{-x}$ have linear resolution, and $\reg(I:x) = d-1$, $\reg(I_{-x})=d$. By \Cref{exact}, 
\begin{gather*}
    \reg(I) \leq \max\{ \reg (I : x) + 1, \reg (I,x)\} = \max\{ d, \reg (I,x)\}. 
\end{gather*}

Now, $(I,x) = (I_{-x},x)$, so $\reg(I,x)= \reg(I_{-x},x)$. Further, $\reg(I_{-x},x) \le \reg(I_{-x})$ by proof of \cite[Lemma 2.10]{DHS13}. Hence, $\reg(I,x) \le \reg(I_{-x}) = d$. Thus, $\reg(I)\le d$.
If $I_{-x}=\la 0 \ra$, then $I=x (I:x)$ and since $(I:x)$ has linear resolution by induction, the claim follows by \Cref{lem_prodlinear}.
 \end{proof}

Next, we give a proof for the Betti splitting of an equigenerated good-vertex splittable ideal by applying the following result by Christopher Francisco, Huy T\`ai H\`a and Adam Van Tuyl. 

\begin{lemma} \cite[Corollary 2.4]{FHVT09} 
Let $I$ be a monomial ideal in $S$, and suppose that $J$ and $K$ are monomial
ideals in $S$ such that $\G(I)$ is the disjoint union of $\G(J)$ and $\G(K)$. If both $J$ and $K$ have linear
resolution, then $I = J + K$ is a Betti splitting.
\end{lemma}

\begin{proposition} \label{prop_gvsbetti}
   Let $I$ be an equigenerated ideal such that $I= x\cdot(I:x) + I_{-x}$ is a good-vertex splitting. Then this is also a Betti splitting.
\end{proposition}

\begin{proof}
From the definition, it is clear that $\G(I)$ is the disjoint union of $\G(x(I:x))$ and $\G(I_{-x})$. By Theorem \ref{gvs_linres}, both $(I:x)$ and $I_{-x}$ have linear resolution and, $x(I:x)$ has linear resolution by \Cref{lem_prodlinear}. The conclusion now follows by the above lemma.
\end{proof}

\section{Families of good-vertex splittable  ideals}

In this section, we explore classes of good-vertex splittable ideals. Since an equigenerated good-vertex splittable ideal has linear resolution, this gives us an easy combinatorial tool to establish linear resolution without explicitly computing the syzygies. 

\begin{example} \label{ex_borel}
     A monomial ideal $I$ is said to be \emph{Borel} if it satisfies the condition that for every monomial $\fu \in I$ and every $x_j$ dividing $\fu$, $\dfrac{\fu}{x_j}\cdot x_i \in I$ whenever $i<j$. Shalom Eliahou and Michel Kervaire proved that an equigenerated Borel ideal has linear resolution by explicitly constructing a resolution in \cite{EK90}. For a Borel ideal $I \subset \s$, $x_1$ is a good vertex, and  $I=x_1 (I:x_1) ~+~I_{-x_1}$ is a good-vertex splitting. This was shown by Marilena Crupi and Antonino Ficarra in \cite[Proposition 1(a)]{CF24}.
\end{example}

\begin{remark}\label{Rem_doubleinduction}
    To show good-vertex splittability, one can often apply induction on the lexicographical order on the pair $(d,\ell)$, where $d= \maxd(I)$ and $\ell$ is the cardinality of $\supp(I)$, for a monomial ideal $I$. Note that the cases $d=0$ and $\ell=0$ are trivially good-vertex splittable. The case $d=1$ was shown in \Cref{gvs_ex} (1), and the case $\ell=1$ is a principal ideal which is good-vertex splittable. Further, if $I$ is squarefree and $\ell=2$, then $I$ is generated by variables or $I= \la x_1x_2 \ra$, making it good-vertex splittable.
\end{remark}

\begin{remark} \label{rmk_redsqf}
    \noindent For $I$ a squarefree monomial ideal generated in degree $d$,
\begin{enumerate}
    \item Consider the combinatorial visualization discussed before \Cref{ex_goodvertex}.\\ Let $\supp(I)=\ell$. If $\Delta$ is taken as the $(d-1)$-skeleton of the simplex on $\ell$-vertices, the facets corresponding to monomials in $\G(I)$ are colored blue, and the rest in red. The red facets correspond to monomials in $\G(\cp(I))$.
    \item  Note that for any good vertex $x \in \supp(I)$, $\supp((I:x)) \subseteq \supp(I) \setminus \{x\}$. Any monomial $\fq$ is a generator of $\cp((I:x))$ if and only if $\deg(\fq)=d-1$, $\supp(\fq) \subset \supp((I:x))$ and, $\fq \cdot x \notin I$, i.e., $\fq \cdot x$ corresponds to a red facet of $\Delta$.
\end{enumerate}  
\end{remark}
  
In \cite{DNSVTW24}, Anton Dochtermann, 
 Ritika Nair, Jay Schweig, Adam Van Tuyl and Russ Woodroofe prove that pure simplicial complexes having sufficiently many facets are vertex decomposable. \Cref{dual_vd} is the Alexander dual version of this result \cite[Theorem 1]{DNSVTW24}. The proof presented is an adaptation of \cite{DNSVTW24} and gives a nice application of the combinatorial visualization discussed above. Our proof relies on simple combinatorial techniques. This also recovers \Cref{thm_TY}, a result first proved by Naoki Terai and Ken-ichi Yoshida using homological methods.

\begin{theorem} \label{dual_vd}
Let $I$ be a squarefree monomial ideal generated in degree $d \ge 1$ such that $|\G(\cp(I))|\leq 2d-1$. Then $I$ is good-vertex splittable.
\end{theorem}

\begin{proof}
 For a squarefree monomial ideal $I$ generated in degree $d$ and having support of size $\ell$, \Cref{Rem_doubleinduction} establishes good-vertex splittability for $d \le 1$ or $\ell \le 2$. We apply induction on $(d,\ell)$. For a fixed $(d,\ell)$ with $d \ge 2$ and $\ell \ge 3$, and for all $(a,b)<(d,\ell)$, suppose $I$ is good-vertex splittable whenever $I$ is a squarefree monomial ideal generated in degree $a$, and having support of size $b$ such that $|\G(\cp(I))| \le 2a-1$. Let $I$ be a squarefree monomial ideal generated in degree $d$ and having $|\supp(I)|=\ell$, with $|\G(\cp(I))| \le 2d-1$. We need to show the existence of a good vertex $x$ such that\\
(i) $(I:x)$ generated in degree $d-1$ satisfies $|\G(\cp((I:x)))|\le 2(d-1)-1=2d-3$ and \\
(ii) $I_{-x}$ satisfies $|\G(\cp(I_{-x}))|\le 2d-1$. This is automatically true for any $x \in \supp(I)$ since $\G(\cp(I_{-x}))$ in $\k[\hat{x}]$ is a subset of  $\G(\cp(I))$ in $\s$. \\

We are left to establish the existence of a good vertex and (i) above.
We use the geometric realization from \Cref{rmk_redsqf}. Let $\Delta$ be the $(d-1)$-skeleton of the simplex on $\ell$ vertices. For each bad pair $(y,\fu)$, $\D_{(y,\fu)}$ is the $(d-1)$-skeleton of the simplex on the set $\supp(\fu) \cup \{y\}$ of $d+1$ vertices with all the facets containing $y$ colored red, and the facet corresponding to $\fu$ colored blue, i.e., there are $d$ red facets and $1$ blue facet for $\D_{(y,\fu)}$. The following cases are possible:\\

(a) There are at least two bad vertices for $I$, say, $y,z$. Let $(y,\fu)$ and $(z,\fv)$ be two bad pairs with the corresponding complexes $\Delta_{(y,\fu)}$ and $\Delta_{(z,\fv)}$ respectively. Each of these complexes has $d$ red facets. Since the number of red facets of $\D$ is at most $2d-1$, $\Delta_{(y,\fu)}$ and $\Delta_{(z,\fv)}$ must have common red facets and by \Cref{prop_commonred}, there is exactly one common red facet corresponding to a monomial of the form $ \frac{\fu}{x}\cdot y$ for some $x \in \supp(\fu)$. The total number of red facets from $\Delta_{(y,\fu)}$ and $\Delta_{(z,\fv)}$ is $2d-1$, so $\D$ has no more red facets as $|\G(\cp(I))|\leq 2d-1$. Now observe that the vertex $x$ appears in $d-1$ red facets of $\Delta_{(y,\fu)}$ and is not present in $\Delta_{(z,\fv)}$ by the proof of \Cref{prop_commonred}. So, $x$ appears in exactly $d-1$ red facets of $\D$, and hence is a good vertex since for any monomial $\fw \in \G(I)$, $(x, \fw)$ is a bad pair if and only if $\Delta_{(x,\fw)}$ has $d$ red facets containing $x$. Let $W:=\supp((I:x)) \subseteq \supp(I) \setminus \{x\}$. By Part (2) of \Cref{rmk_redsqf}, any monomial $\fq \in \G(\cp((I:x)))$ if and only if $\fq \cdot x$ gives a red facet in $\Delta$ and this happens if and only if $\fq \cdot x$ gives a red facet in $\Delta_{(y,\fu)}$. Since there are exactly $d-1$ red facets in $\Delta_{(y,\fu)}$ that contain $x$, $|\G(\cp((I:x)))| \le d-1\leq 2d-3$, as required.\\

(b) There is exactly one bad vertex for $I$, say $y$. Let $(y,\fu)$ be a bad pair. Since $I$ is squarefree, $\ell \ge d$. If $\ell=d+1$, all the vertices are in $\D_{(y,\fu)}$ with $d$ red facets containing $y$ and one blue facet corresponding to $\fu$. But this means $y \notin \supp(I)$ and $I$ is the principal ideal $\la \fu\ra$. So, let $\ell \geq d+2$. There are $\ell-1$ good vertices. Suppose for every good vertex $x$, $|\G(\cp((I:x)))|> 2d-3$. Then, by Part (2) of \Cref{rmk_redsqf}, each good vertex is present in at least $2d-2$ red facets. We have $y$ in $d$ red facets. Let us count the total occurrence of vertices in red facets.  $|\G(\cp(I))|\leq 2d-1$ and each red facet has $d$ vertices, so,
\[(2d-1)d \geq (\ell-1)(2d-2)+d \geq (d+1)(2d-2)+d=2d^2+d-2\]   
which is not possible for $d>1$.\\

(c) There are no bad vertices. If $\ell=d$, $I$ is a principal ideal and we are done. So, let $\ell\geq d+1$. If there is some vertex $x$ present in at most $2d-3$ red facets, this gives $|\G(\cp((I:x)))| \le 2d-3$, and we are done. Suppose each (good) vertex appears in at least $2d-2$ red facets, the same counting argument as above gives
\[(2d-1)d \geq \ell(2d-2)\geq (d+1)(2d-2)=2d^2-2\]
This leaves us to verify the case forcing equality throughout giving $d=2$ and $\ell=3$. There are two possibilities, $I=\la x_1x_2,x_1x_3\ra$ and $I=\la x_1x_2,x_1x_3,x_2x_3\ra$, both of which are good-vertex splittable. 
\end{proof}

 Let $I$ be a squarefree monomial ideal generated in degree $d$ and $V \subseteq \supp(I)$. The complement of $I$ with respect to $V$, $\cp_V(I):= \la \fu \in \G(\cp(I))~|~ \supp(\fu) \subset V \ra$. Simply put, $\cp_V(I)$ is generated by all the squarefree degree $d$ monomials supported on $V$ that are not in $I$. For example, for $I=\la x_1x_2,x_2x_3,x_3x_4\ra$ as in \Cref{ex:basic}, let $V=\{x_1,x_3,x_4\}$. Then $\cp_V(I)=\la x_1x_3,x_1x_4 \ra$. 

\begin{theorem}\label{complement2}
Let $I$ be a squarefree monomial ideal generated in degree $d$ such that for each $V \subset \supp(I)$ with $|V| \ge d$, we have $|\G(\cp_V(I))| \leq \dfrac{3}{(d+1)}|V|-1$. Then $I$ is good-vertex splittable. 
\end{theorem}

\begin{proof}
   Let $|\supp(I)| > 2$ and $d \ge 2$. Let $\F$ be a family of equigenerated ideals such that for each $V\subset \supp(I)$ with $|V| \ge d$, we have $|\G(\cp_V(I))|\leq \dfrac{3}{(d+1)}|V|-1$. If $d=2$, $I$ is the edge ideal of a graph $G$ on the set of vertices given by $\supp(I)$. The condition says that the complement of any induced subgraph of $G$ on a subset $V$ contains at most $\dis |V|-1$ edges and so is chordal. Hence, $I$ is good-vertex splittable  by \Cref{thm_chordal_vs}, \Cref{vs_implies_gvs}. So, we assume $d\geq 3$. By \Cref{defn_gvs_family}, we need to show the existence of a good vertex $x$ such that $(I:x) \in \F$ and $I_{-x} \in \F$. The second condition holds for any $x$ since $I_{-x}$ is generated in degree $d$ and $V \subset \supp(I_{-x}) \implies V \subset \supp(I)$.\\
   
Suppose $y$ is a bad vertex for $I$ generated in degree $d \ge 3$. For a bad pair $(y,\fu)$ of $I$, set $V=\supp(\fu)\cup\{y\}$. The complex $\D_{(y,\fu)}$ on $V$ has $d$ red facets and one blue facet, i.e.,  $|\G(\cp_V(I))| = d$. By the definition of $I \in \F$, $\dis |\G(\cp_V(I))| \le \frac{3}{(d+1)}|V|-1 = \frac{3}{(d+1)}(d+1)-1=2$, i.e., $d \le 2$, a contradiction. So, every vertex is a good vertex.\\

Choose a good vertex $x$. Let $J:=(I:x)$, generated in degree $d-1$ by \Cref{prop_good}. Let $W \subset \supp(J)$ with $|W| \ge d-1$. We need to show that $ |\G(\cp_W(J))| \leq \frac{3}{d}|W|-1$. If $|W|=d-1$, then $|\G(\cp_W(J))|$ is at most one, and the above inequality holds as $d>2$ here. Suppose $|W| \ge d$. $\G(\cp_W(J))$ consists of some monomials $\fu$ in degree $d-1$ such that $\fu \cdot x \notin I$, i.e., $\fu \cdot x \in \G(\cp(I))$. Hence,
\[
   |\G(\cp_W(J))| \le  |\G(\cp_{W \cup \{x\}}(I))| 
    \le   \frac{3}{(d+1)}(|W|+1)-1 \text{, since  } I \in \F
    \]
and
\[ 
\frac{3}{(d+1)} (|W|+1) -1\le \frac{3}{d}|W|-1 \iff d \le |W|
\]
which is true.
Thus, $J = (I:x) \in \F$. 
 \end{proof}

\section{Edge ideals with regularity at most $3$}

In this section, we first make the observation that a graph having too many good vertices has edge ideal with regularity $2$.

 \begin{theorem} \label{gvreg}
    
     If the edge ideal of a graph has at most $3$ bad vertices, it has a linear resolution.

\end{theorem}

\begin{proof}
    Let $I=I(G)$ be the edge ideal of a graph $G$. By Fr\"oberg's Theorem (\Cref{thm_chordal}), $I$ does not have a linear resolution if the graph complement $G^\c$ has an induced cycle $C_n$ on $n \ge 4$ vertices. Then these corresponding $n$ vertices are bad for $I$ by \Cref{ex_goodvertex}(2) resulting it at least $4$ bad vertices for $I$.
\end{proof}

In the remainder of this section, we present some results for edge ideals having regularity at most $3$ by giving sufficient conditions on induced subgraphs of size $5$ or $6$. Before proceeding, we quote the following lemma which is a version of \Cref{exact} in the case of edge ideals.
\begin{lemma} \cite[Lemma 3.1]{DHS13}\label{lem_reg_edge}
    Let $G$ be a graph and $x$ be a vertex of $G$ with $\ngb_G(x) \ne \emptyset$. Then, 
    \[
\reg(I(G)) \leq \max\{ \reg (I(G \setminus \st_G(x))) + 1, \reg (I(G \setminus \{x\}))\}. 
\]
Moreover,  $\reg(I(G))$ is {\it equal} to one of these terms.
\end{lemma}

\begin{theorem}
\label{reg3graphfive}
    Let $G$ be a simple graph without isolated vertices such that the induced subgraph of $G$ on any five vertices  is neither $K_1 \sqcup 2K_2$ nor  $K_1 \sqcup P_4$.  Then the edge ideal of $G$ has regularity at most $3$. 
\end{theorem}

\begin{proof}
     We prove the theorem by induction on $n$, the number of vertices of $G$. It is known that if $n\leq 5$, the regularity of $I(G)$ is at most 3. Assume that the theorem is true up to some $n\geq 5$. Let $G$ be a graph with $n+1$ vertices satisfying the hypotheses. Let $x$ be a vertex of $G$. We use the above lemma. Note that the graph $G \setminus \{x\}$ has $n$ vertices and also satisfies the hypotheses of the theorem, so by induction, $\reg (I(G \setminus \{x\})) \le 3$. It suffices to show that $\reg (I({G \setminus \st_G(x)})) \leq 2$, i.e., $(G \setminus \st_G(x))^{\c}$ is chordal by \Cref{thm_chordal}. Suppose on the contrary, the graph $(G \setminus \st_G(x))^{\c}$ has an induced cycle on the vertices $x_1,x_2, \ldots, x_t$, where $t\geq 4$. If $t=4$, $\{x_1,x_3\}$ and $\{x_2,x_4\}$ are edges of $G$. This makes the induced subgraph on the five vertices $x,x_1,x_2,x_3,x_4$ a  $K_1 \sqcup 2K_2$ in $G$, a contradiction. \\
    Next, let $t\geq 5$. Then, $\{x_1,x_3\}$, $\{x_1,x_4\}$, $\{x_2,x_4\}$ are edges of $G$, so the induced subgraph of $G$ on the five vertices $x,x_1,x_2,x_3,x_4$ is $K_1 \sqcup P_4$, again a contradiction.
\end{proof}

\begin{corollary}
    Let $G$ be a simple graph without isolated vertices. If $G$ is gap-free and $P_4$-free, then the edge ideal of $G$ has regularity at most 3.  
\end{corollary}

\begin{theorem}
\label{reg3-graph-six}
    Let $G$ be a simple graph with no isolated vertices. Suppose that the  induced subgraph of $G$ on any six vertices  is not one of the graphs in Table \ref{tab:six}.   Then the edge ideal $I(G)$ of $G$ has regularity at most 3.

    \begin{table}[h!]
        \centering
        \begin{tabular}{|m{3cm}|m{3cm}|m{3cm}|m{3cm}|}
        \hline 
              \begin{tikzpicture}[scale=0.7]
            
                \tkzDefPoint(0,0){1}
                \tkzDefPoint(0,1){2}
                
                \tkzDefPoint(1,0){3}                
                \tkzDefPoint(1,1){4}
                
                \tkzDefPoint(2,0){5}
                \tkzDefPoint(2,1){6} 

                \draw [thick] (1) -- (2) ; 
                \draw [thick] (3) -- (4); 
                \draw [thick] (5) -- (6);
                
                \tkzDrawPoints[size=3](1,2,3,4,5,6)

                %\node[below] at (2,0) {{\normalfont (a)}};
            \end{tikzpicture} & 
            \begin{tikzpicture}[scale=0.7]
            
                \tkzDefPoint(0,0){1}
                \tkzDefPoint(1,0){2}
                \tkzDefPoint(2,0){3}
                \tkzDefPoint(3,0){4}
                
                \tkzDefPoint(1,1){5}                
                \tkzDefPoint(2,1){6}

                \draw [thick] (1) -- (4) ;  
                \draw [thick] (5) -- (6);
                
                \tkzDrawPoints[size=3](1,2,3,4,5,6)

                %\node[below] at (2,0) {{\normalfont (a)}};
            \end{tikzpicture} &
            \begin{tikzpicture}[scale=0.7]
            
                \tkzDefPoint(0,0){1}
                \tkzDefPoint(1,0){2}
                \tkzDefPoint(2,0.5){3}
                \tkzDefPoint(2,-0.5){4}
                
                \tkzDefPoint(3,0.5){5}                
                \tkzDefPoint(3,-0.5){6}

                \draw [thick] (1) -- (2) -- (3) -- (4) -- (2) ;  
                \draw [thick] (5) -- (6);
                
                \tkzDrawPoints[size=3](1,2,3,4,5,6)

                %\node[below] at (2,0) {{\normalfont (a)}};
            \end{tikzpicture} & 
            \begin{tikzpicture}[scale=0.7]
            
                \tkzDefPoint(0,0){1}
                \tkzDefPoint(1,0){2}
                \tkzDefPoint(2,0){3}
                \tkzDefPoint(3,0){4}
                
                \tkzDefPoint(1,1){5}                
                \tkzDefPoint(1,2){6}

                \draw [thick] (1) -- (4) ;  
                \draw [thick] (2) -- (6);
                
                \tkzDrawPoints[size=3](1,2,3,4,5,6)

                \node[above] at (1,2.2) {}; %fake node

                %\node[below] at (2,0) {{\normalfont (a)}};
            \end{tikzpicture}\\ [0.5ex]
            \hline 
              \begin{tikzpicture}[scale=0.7]
            
                \tkzDefPoint(0,0){1}
                \tkzDefPoint(1,0){2}
                \tkzDefPoint(2,0){3}
                \tkzDefPoint(3,0){4}
                
                \tkzDefPoint(0.5,1){5}                
                \tkzDefPoint(1.5,1){6}

                \draw [thick] (1) -- (4) ;  
                \draw [thick] (2) -- (5) -- (6) -- (2);
                
                \tkzDrawPoints[size=3](1,2,3,4,5,6)

                %\node[below] at (2,0) {{\normalfont (a)}};
            \end{tikzpicture} &
            \begin{tikzpicture}[scale=0.35]
            
                \tkzDefPoint(0,0){1}
                \tkzDefPoint(0,2){2}
                \tkzDefPoint(2,3){3}
                
                \tkzDefPoint(4,2){4}
                \tkzDefPoint(4,0){5}
                \tkzDefPoint(2,1){6} 

                \draw [thick] (1) -- (2) -- (5) -- (4) -- (1); 
                \draw [thick] (3) -- (6); 
                
                \tkzDrawPoints[size=3](1,2,3,4,5,6)

                \node[above] at (2,3.2) {};

                %\node[below] at (2,0) {{\normalfont (a)}};
            \end{tikzpicture}
            & 
            \begin{tikzpicture}[scale=0.35]
            
                \tkzDefPoint(0,0){1}
                \tkzDefPoint(1,0){2}
                \tkzDefPoint(2,0){3}                
                \tkzDefPoint(3,0){4}
                \tkzDefPoint(4,0){5}
                \tkzDefPoint(5,0){6} 

                \draw [thick] (1) -- (6); 
                
                \tkzDrawPoints[size=3](1,2,3,4,5,6)
                
                %\node[below] at (2,0) {{\normalfont (a)}};
            \end{tikzpicture}
               &   
               \begin{tikzpicture}[scale=0.7]
            
                \tkzDefPoint(0,0){1}
                \tkzDefPoint(1,0){2}
                \tkzDefPoint(2,0.5){3}                
                \tkzDefPoint(1,1){4}
                \tkzDefPoint(0,1){5}
                \tkzDefPoint(3,0.5){6} 

                \draw [thick] (1) -- (2) -- (3) -- (4) -- (5) -- (1); 
                \draw [thick] (3) -- (6); 
                
                \tkzDrawPoints[size=3](1,2,3,4,5,6)

                %\node[below] at (2,0) {{\normalfont (a)}};
            \end{tikzpicture} \\ [0.5ex]
               \hline  
              \begin{tikzpicture}[scale=0.7]
            
                \tkzDefPoint(0,0){1}
                \tkzDefPoint(0,1){2}
                
                \tkzDefPoint(1,0){3}                
                \tkzDefPoint(1,1){4}
                
                \tkzDefPoint(2,0){5}
                \tkzDefPoint(2,1){6} 

                \draw [thick] (1) -- (2) -- (4) -- (3) -- (1) ; 
                \draw [thick] (3) -- (5);  
                \draw [thick] (4) -- (6); 
                
                \tkzDrawPoints[size=3](1,2,3,4,5,6)

                \node[above] at (2,1.1) {};
            \end{tikzpicture} &
            \begin{tikzpicture}[scale=0.7]
            
                \tkzDefPoint(0,0){1}
                \tkzDefPoint(1,0){2}
                
                \tkzDefPoint(2,0.5){3}                
                \tkzDefPoint(2,-0.5){4}
                
                \tkzDefPoint(3,0.5){5}
                \tkzDefPoint(3,-0.5){6} 

                \draw [thick] (1) -- (2) -- (3) -- (4) -- (2) ; 
                \draw [thick] (3) -- (5);  
                \draw [thick] (4) -- (6); 
                
                \tkzDrawPoints[size=3](1,2,3,4,5,6)

                %\node[below] at (2,0) {{\normalfont (a)}};
            \end{tikzpicture}&
            \begin{tikzpicture}[scale=0.7]             
                \tkzDefPoint(0,0){1}
                \tkzDefPoint(1,0){2}
                
                \tkzDefPoint(2,0.5){3}                
                \tkzDefPoint(1,1){4}
                
                \tkzDefPoint(0,1){5}
                \tkzDefPoint(3,0.5){6} 

                \draw [thick] (1) -- (2) -- (3) -- (4) -- (5) -- (1) ; 
                \draw [thick] (1) -- (3) -- (6);   
                
                \tkzDrawPoints[size=3](1,2,3,4,5,6)

                %\node[below] at (2,0) {{\normalfont (a)}};
            \end{tikzpicture} &  
            \begin{tikzpicture}[scale=0.7]
            
                \tkzDefPoint(0,0){1}
                \tkzDefPoint(0,1){2}
                \tkzDefPoint(1,1){3}
                \tkzDefPoint(1,0){4}
                
                \tkzDefPoint(2,0){5}                
                \tkzDefPoint(2,1){6}

                \draw [thick] (1) -- (3) -- (4);  
                \draw [thick] (5) -- (1) -- (2) -- (6);
                
                \tkzDrawPoints[size=3](1,2,3,4,5,6)

                \node[above] at (2,1.1) {};

                %\node[below] at (2,0) {{\normalfont (a)}};
            \end{tikzpicture}
           
            \\ [0.5ex]
        \hline 
            \begin{tikzpicture}[scale=0.7]             
                \tkzDefPoint(0,-0.5){x}
                \tkzDefPoint(0,0){1}
                \tkzDefPoint(1,0){2}
                
                \tkzDefPoint(2,0.5){3}                
                \tkzDefPoint(1,1){4}
                
                \tkzDefPoint(0,1){5}
                \tkzDefPoint(3,0.5){6} 

                \draw [thick] (1) -- (2) -- (3) -- (4) -- (5) -- (1) ; 
                \draw [thick] (1) -- (3) -- (6);  
                \draw [thick] (3) -- (5);
                
                \tkzDrawPoints[size=3](1,2,3,4,5,6)

                \node[above] at (2,1.1) {};

                %\node[below] at (2,0) {{\normalfont (a)}};
            \end{tikzpicture} & &  &\\ [0.5ex]
            \hline 
        \end{tabular} 
        \caption{}
        \label{tab:six}
    \end{table}
\end{theorem} 

\begin{proof}
    The proof of the theorem is similar to that of \Cref{reg3graphfive} and we prove by induction on $n$, the number of vertices of $G$. We apply \Cref{lem_reg_edge}. It is known that if $n\leq 5$, the regularity of $I(G)$ is at most 3. Assume that the theorem is true up to some $n\geq 5$. Let $G$ be a graph with $n+1$ vertices satisfying the hypotheses. Let $x$ be a vertex of $G$ and $y \in \ngb_G(x)$. It suffices to show that $\reg (I(G \setminus \st_G(x))) \leq 2$, i.e., $(G \setminus \st_G(x))^{\c}$ has no induced cycle of length at least 4. Suppose on the contrary that $x_1,\ldots,x_t$ are vertices in $(G \setminus \st_G(x))^{\c}$ such that the induced subgraph of $G^\c$ on these vertices is a cycle, where $t\geq 4$. We consider two cases. If $t=4$, this gives the induced subgraph of $G$ on the vertices $x,y,x_1,x_2,x_3,x_4$ consisting of the edges $\{x_1,x_3\}$, $\{x_2,x_4\}$, and $\{x,y\}$, no edges between $x$ and the set $x_1, x_2, x_3, x_4$, and possible edges from $y$ to any of the vertices $x_1, x_2, x_3, x_4$. Next, let $t\geq 5$. Then the induced subgraph of $G$ on the vertices $x,y,x_1,x_2,x_3,x_4$ consists of the edges $\{x_1,x_3\}$, $\{x_1,x_4\}$, $\{x_2,x_4\}$, and $\{x,y\}$, no edges between $x$ and the set $x_1, x_2, x_3, x_4$, and possible edges from $y$ to any of the vertices $x_1, x_2, x_3, x_4$. The induced subgraph of $G$ on $x,y,x_1,x_2,x_3,x_4$ from both the cases are listed in \Cref{tab:six}. 
\end{proof}

Most recently, Kanno \cite{Ka26} gave a characterization for the edge ideal of a graph to have regularity at most $3$ in terms of triangulation of closed $2$-manifolds. Our result explicitly gives a class of graphs satisfying Kanno's topological condition: If a graph $G$ does not contain the aforementioned forbidden induced subgraphs on any five (or six) vertices, then the graph complement of $G$ has no subgraph that is a triangulation of an orientable closed $2$-manifold.

\section{Further Questions}

We conclude by mentioning some questions related to the discussion in this paper.\\

This paper introduces the concept of a good vertex. Since a good vertex of a squarefree monomial ideal is the shedding vertex in its Alexander dual complex, this helps translate combinatorial consequences of shedding vertices in simplicial complexes to algebraic terms for squarefree monomial ideals.

\begin{question}
    For a given characteristic, it is not easy to find Cohen-Macaulay complexes which are not vertex decomposable. In the Alexander dual language, we are looking for squarefree monomial ideals having a linear resolution but are not vertex splittable. An interesting direction to pursue would be the use of algebraic techniques to find monomial ideals that have linear resolution but do not satisfy the recursive conditions in \Cref{defn_gvs_family}. But we are not sure how feasible this process is.
\end{question}

\begin{question}

This question is inspired by a communication of the second author with Naoki Terai and Ken-ichi Yoshida. In \cite{TY06}, Naoki Terai and Ken-ichi Yoshida proved the following result using reduced homology. This is now obtained as a consequence of \Cref{gvs_linres} using the simple combinatorial proof given for \Cref{dual_vd}.
\begin{corollary}\cite[Theorem 3.3]{TY06} \label{thm_TY}
Let $I$ be a squarefree monomial ideal generated in degree $d$ such that $|\G(\cp(I))|\leq 2d-1$. Then $I$ has linear resolution. 
\end{corollary}

 \Cref{thm_TY} tells us that an equigenerated monomial ideal whose complement ideal has a limited number of generators has a linear resolution. The bound here is sharp (see \cite[Example 13]{DNSVTW24}, \cite[Example 5.2]{TY06}). What happens if we drop the assumption of being equigenerated? By naturally extending \Cref{def_comp} to complement for any squarefree monomial ideal, can we say that a squarefree monomial ideal whose complement ideal has a small enough number of generators is componentwise linear? 
 
\end{question}

\begin{question}
    \Cref{gvreg} establishes that if for the edge ideal of a graph, almost all vertices are good, then the regularity is $2$. If the edge ideal has $t \ge 4$ many bad vertices, what can we say about the regularity of the ideal? Also, what is the scope of extending \Cref{gvreg} to ideals generated in degree $> 2$?
\end{question}

\bibliographystyle{amsplain}
\bibliography{bibliography}

\end{document}